\documentclass{article}

\usepackage{amsmath,amssymb}
\usepackage[dvips]{epsfig}
\usepackage[latin1]{inputenc}

\newcommand{\ZZ}{\mathbb{Z}}
\newcommand{\CC}{\mathbb{C}}
\newcommand{\QQ}{\mathbb{Q}}
\newcommand{\Det}{\operatorname{Det}}

\newcommand{\magm}{\operatorname{Gri}}
\newcommand{\bess}{\operatorname{Bess}}
\newcommand{\rt}{\operatorname{root}}

\newcommand{\LLL}{\operatorname{LLL}}
\newcommand{\YY}{\operatorname{Y_{Y} Y}}
\newcommand{\supp}{\operatorname{Supp}}

\newtheorem{theo}{Theorem}[section]
\newtheorem{prop}[theo]{Proposition}
\newtheorem{lemma}[theo]{Lemma}

\newtheorem{question}{Question}

\newenvironment{proof}{\begin{trivlist}\item{\bf{Proof.}}}
  {\hfill\rule{2mm}{2mm}\end{trivlist}}

\title{A Hopf operad of forests of binary trees and related
  finite-dimensional algebras}
\date{\today}
\author{Frédéric Chapoton}

\begin{document}

\maketitle

\begin{abstract}
  The structure of a Hopf operad is defined on the vector spaces
  spanned by forests of leaf-labeled, rooted, binary trees. An explicit
  formula for the coproduct and its dual product is given, using a
  poset on forests.
\end{abstract}

\setcounter{section}{-1}

\section{Introduction}

The theme of this paper is the algebraic combinatorics of leaf-labeled
rooted binary trees and forests of such trees. We shall endow these
objects with several algebraic structures.

The main structure is an operad, called the Bessel operad, which is
the suspension of an operad defined by a distributive law between the
suspended commutative operad and the operad of commutative
non-associative algebras (sometimes called Griess algebras). The
Bessel operad may be seen as an analog of the Gerstenhaber operad
\cite{voronov}, which is the suspension of an operad defined by a
distributive law between the suspended commutative operad and the
Poisson operad. Unlike the Gerstenhaber operad, the Bessel operad has
a simple combinatorial basis, given explicitly by forests of
leaf-labeled rooted binary trees.

The Bessel operad, like the Gerstenhaber operad, is a Hopf operad.
More precisely, they are both endowed with a cocommutative coproduct.
This gives rise to a family of finite-dimensional coalgebras. In the
dual vector spaces of the Bessel operad, one gets algebras based on
forests of leaf-labeled binary trees.

An explicit formula is obtained for the coproduct in these coalgebras
of forests (and therefore for their dual products), using a poset
structure on the set of forests, which may be of independent interest.

The first section is devoted to the definition of a distributive law
between the suspended commutative operad and the Griess operad. The
suspension of the operad defined by this distributive law is
introduced in the next section. The coproduct is defined and shown to
be given by an explicit sum in the third section. In the fourth
section, the dual algebras are briefly studied.

\section{A distributive law}

All the operads considered here are in the monoidal category of
complexes of vector spaces over $\QQ$ with zero differential,
\textit{i.e.} the category of vector spaces over $\QQ$ which are
graded by $\ZZ$, with Koszul sign rules for the tensor product. An
operad \cite{ginzkapr,livreMSS,may} is seen through its underlying
functor from the groupoid of finite sets to this monoidal category. An
Hopf operad is an operad $\mathcal{P}$ with a coassociative morphism
of operad from $ \mathcal{P} $ to $\mathcal{P} \otimes \mathcal{P}$.

\begin{figure}
  \begin{center}
    \leavevmode 
    \epsfig{file=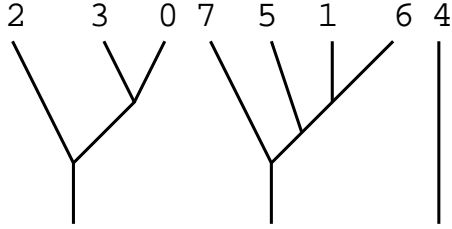,width=6cm} 
    \caption{A forest on $\{0,1,2,\dots,7\}.$}
    \label{exforest}
  \end{center}
\end{figure}

A \textit{tree} is a leaf-labeled rooted binary tree and a
\textit{forest} is a set of such trees, see Fig. \ref{exforest}.
Vertices are either inner vertices (valence $3$) or leaves and roots
(valence $1$). By convention, edges are oriented towards the root.
Leaves are bijectively labeled by a finite set. An half-edge is a pair
made of an inner vertex and an incident edge (incoming or outcoming).
Trees and forests are pictured with their roots down and their leaves
up, but are not to be considered as planar.

\subsection{The determinant operad and orientations}

An \textit{orientation} of a finite set $X$ is a maximal exterior
power of the elements of this set, \textit{i.e} a generator of the
$\ZZ$-module $\Lambda^{|X|} \ZZ X$.

Let us recall the definition of the suspended commutative associative
operad $\Det$ introduced by Ginzburg and Kapranov \cite{ginzkapr}. Let
$I$ be a finite set, then $\Det(I)$ is the determinant vector space of
$\QQ I$. This vector space is one-dimensional, spanned by the
orientations of $I$ (for example $1\wedge 3\wedge 4\wedge 2$ in
$\Det(\{1,2,3,4\})$) and is placed in degree $|I|-1$. The composition
of the operad $\Det$ is given by the rule
\begin{equation}
 (x\wedge \star) \circ_\star y =x\wedge y,
\end{equation}
for all $x\in\Det(I)$ and $y\in\Det(J)$.

It is well known and easy to check that $\Det$ has the presentation by
the antisymmetric generator $e_{i,j}=i \wedge j$ of degree $1$ in
$\Det(\{i,j\})$ satisfying
\begin{equation}
  \label{reldet}
  e_{i,\star}\circ_\star e_{j, k} =e_{k, \star }\circ_\star e_{i ,
    j}.
\end{equation}
The operad $\Det$ is binary quadratic and Koszul, see \cite{ginzkapr}
for the definitions of these notions.

\subsection{The Griess operad and rooted binary trees}

The operad $\magm$ describing commutative but not necessarily
associative algebras (sometimes called Griess algebras) admits the
following description. The space $\magm(I)$ has a basis indexed by
rooted binary trees with leaves labeled by $I$ and the composition is
grafting. This vector space is placed in degree $0$. In fact,
$\magm$ is the free operad on a binary symmetric generator
$\omega_{i,j}$ of degree $0$ corresponding to the unique rooted binary
tree with two leaves labeled by $\{i,j\}$. The operad $\magm$ is
binary quadratic and Koszul.

\subsection{The operad $B$ of root-oriented forests}

For the definition and properties of the notion of distributive law
from an operad to another one, see \cite{markl}.

\begin{prop}
  The following formula defines a distributive law \\ from $\magm
  \circ \Det $ to $\Det \circ \magm$ :
  \begin{equation} \label{distri} \omega_{i,\star} \circ_\star e_{j,k}
    = e_{j,\star} \circ_\star \omega_{i,k} - e_{k,\star} \circ_\star
    \omega_{i,j}.
  \end{equation}\end{prop}
\begin{proof}
  As $\magm$ is a free operad, one has only to check that the
  rewriting of 
  \begin{equation}
    \label{differ}    \omega_{i,\star} \circ_\star (e_{j,\#} \circ_{\#}
    e_{k,\ell})    -\omega_{i,\star} \circ_\star (e_{k,\#} \circ_{\#}
    e_{\ell,j}),
  \end{equation}
  using (\ref{distri}) as a replacement rule, gives zero modulo the
  relation (\ref{reldet}) which defines $\Det$. Indeed, one has
  \begin{align*}
    \omega_{i,\star}  & \circ_\star (e_{j,\#} \circ_{\#} e_{k,\ell})=
    (\omega_{i,\star} \circ_\star e_{j,\#}) \circ_{\#} e_{k,\ell}\\&= (
    e_{j,\star} \circ_\star \omega_{i,\#} - e_{\#,\star} \circ_\star
    \omega_{i,j} )\circ_{\#} e_{k,\ell}\\&= e_{j,\star} \circ_\star
    (\omega_{i,\#} \circ_{\#} e_{k,\ell}) -(e_{\#,\star} \circ_\star
    \omega_{i,j} )\circ_{\#} e_{k,\ell}\\&= e_{j,\star} \circ_\star (
    e_{k,\#} \circ_\# \omega_{i,\ell} - e_{\ell,\#} \circ_\#
    \omega_{i,k} ) -(e_{\#,\star} \circ_{\#} e_{k,\ell})\circ_\star
    \omega_{i,j} \\&= ( e_{j,\star} \circ_\star e_{k,\#}) \circ_\#
    \omega_{i,\ell} - (e_{j,\star} \circ_\star e_{\ell,\#} )\circ_\#
    \omega_{i,k} -(e_{\#,\star} \circ_{\#} e_{k,\ell})\circ_\star
    \omega_{i,j} \\&= ( e_{j,\star} \circ_\star e_{k,\#}) \circ_\#
    \omega_{i,\ell} + (e_{j,\star} \circ_\star e_{\#,\ell} )\circ_\#
    \omega_{i,k} +(e_{\#,\star} \circ_{\star} e_{k,\ell})\circ_\#
    \omega_{i,j} \\&= ( e_{j,\star} \circ_\star e_{k,\#}) \circ_\#
    \omega_{i,\ell} + (e_{\ell,\star} \circ_\star e_{j,\#} )\circ_\#
    \omega_{i,k} +(e_{k,\star} \circ_{\star} e_{\ell,\#})\circ_\#
    \omega_{i,j}.
\end{align*}
This expression is invariant by cyclic permutation of $j,k,\ell$.
This shows that the rewriting of (\ref{differ}) is zero, which proves
the proposition.
\end{proof}

Let us summarize the description of the operad defined by this
distributive law.

\begin{prop}
  The operad $B$ defined on $\Det\circ \magm$ by this distributive law
  is isomorphic to the quotient of the free operad generated by
  $e_{i,j}$ antisymmetric in degree $1$ and $\omega_{i,j}$ symmetric
  in degree $0$ by the following relations.
  \begin{align}
    e_{i,\star}\circ_\star e_{j, k} &=e_{k, \star }\circ_\star e_{i ,
      j}, \\
    \omega_{i,\star} \circ_\star e_{j,k} &= e_{j,\star} \circ_\star
    \omega_{i,k} - e_{k,\star} \circ_\star \omega_{i,j}.
  \end{align}
\end{prop}

\begin{figure}  \begin{center}
    \leavevmode 
    \epsfig{file=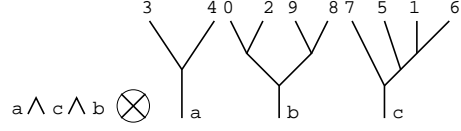,width=6cm} 
    \caption{A
      root-oriented forest on $\{0,1,2,\dots,9\}$.}\label{forexample}
  \end{center}\end{figure}

A root-orientation of a forest $F$ is an orientation of the set of
roots of $F$. A \textit{root-oriented forest} is a tensor product of a
root-orientation and a forest, see Fig. \ref{forexample}. By the
construction of $B$ by a distributive law, the vector space $B(I)$ has
a basis indexed by root-oriented forests. The degree of a
root-oriented forest is the number of roots minus one.

\begin{prop}
  The operad $B$ is binary quadratic and Koszul.
\end{prop}
\begin{proof}
  Koszulness follows from a theorem of Markl \cite{markl} since $B$
  is defined by a distributive law between two Koszul operads.
\end{proof}

Here is a description of the composition in the generators. The
generator $e_{i,j}$ acts on forests by disjoint union. Let $F_1 \sqcup
F_2$ be the disjoint union of two forests $F_1$ and $F_2$. We use
(from now on) the abuse of notation $(-1)^{x}$ for $(-1)^{\deg(x)}$
when $x$ is homogeneous and also $(-1)^{o}$ instead of
$(-1)^{\deg(o)}$ for any kind of orientation $o$. The degree of an
orientation is the number of wedge signs that it contains.

\begin{prop}
\label{action_e}
  Let $o_1 \otimes F_1$ and $o_2 \otimes F_2$ be two root-oriented
  forests. Then
  \begin{equation} e_{\star,\#}\circ_\star (o_1 \otimes
    F_1) \circ_\# (o_2 \otimes F_2)= (-1)^{o_1} o_1 \wedge o_2 \otimes
    (F_1 \sqcup F_2).
  \end{equation}
\end{prop}

\begin{proof}
  The proposition can be restated as follows. Let $x\in B(I)$ and
  $y\in B(J)$. Then
  \begin{equation*}  
    (e_{\star,\#}\circ_\star x) \circ_\# y=(-1)^{x} x \wedge y,  
  \end{equation*}
  Indeed, one has $e_{\#,\star} \circ_\star x=\# \wedge x$ and
  $(x\wedge \#) \circ_\# y=x \wedge y$ by the composition rule of
  $\Det$. The sign is given by $e_{\#,\star}=-e_{\star,\#}$ and
  $\#\wedge x=(-1)^{x+1} x \wedge \#$.\end{proof}

The generator $\omega_{i,j}$ acts on trees by grafting. Let $T_1\vee
T_2$ be the tree obtained by grafting $T_1$ and $T_2$ on the two
leaves of the tree with one inner vertex.

\begin{prop}
  \label{action_o}
  Let $o_1 \otimes T_1$ and $o_2 \otimes T_2$ be two root-oriented
  trees. Then \begin{equation} \omega_{\star,\#}\circ_\star (o_1
    \otimes T_1)\circ_\#(o_2 \otimes T_2) = o \otimes (T_1 \vee T_2),
  \end{equation} where $o$ is the unique root-orientation of the tree
  $T_1 \vee T_2$.
\end{prop}
\begin{proof}
  This is just the composition of $\magm$, restated inside $B$, by
  definition of the composition in an operad defined by a distributive
  law.
\end{proof}

\section{The Bessel operad as a suspension}

This section is devoted to the operad $\bess=\Det \otimes B$ which is
a suspended version of $B$. This suspension is necessary for the
definition of a Hopf operad structure in the next section.

The generating series of the operad $\bess$ has for coefficients the
Bessel polynomials \cite{grosswald,krallfrink}, which are known to
count the forests (sets) of rooted leaf-labeled binary trees, hence
the chosen name.

\subsection{Outer and inner orientations}

By its definition, the vector space $\bess(I)$ has a basis indexed by
tensor products $o_1 \otimes o_2 \otimes F$ where $o_1$ is an
orientation of $I$ and $o_2$ is a root-orientation of the forest $F$.
This tensor product of two orientations is called an \textit{outer
  orientation} of $F$. In this section, an alternative description is
given for this kind of orientation, which will be more convenient
later.

\smallskip

A \textit{global orientation} of a forest $F$ is an orientation of the set
$V(F) \sqcup \{R_F\}$, where $V(F)$ is the set of inner vertices of
$F$ and $R_F$ is an auxiliary element. 

A \textit{local orientation} of a forest $F$ at an inner vertex $v$ is an
orientation of its $3$ incident half-edges (which is of course
equivalent to a cyclic order).

An \textit{inner-oriented forest} is a tensor product
$o\otimes\bigotimes_{v \in V(F)} o_v \otimes F$, where $o$ is a global
orientation of the forest $F$ and the $o_v$ are local orientations of
$F$ at its inner vertices. This will from now on be abridged $o
\otimes F$, where $o$ is a global orientation, the local orientations
being implicit. Notice that the order in the product of the local
orientations do not matter, as they have degree $2$.

\smallskip

One can identify an outer orientation $o_1 \otimes o_2$ with an inner
orientation in the following way.
\begin{enumerate}
\item Consider the exterior product $o_1 \wedge R_F \wedge o_2$ where
  $R_F$ is an auxiliary element.
\item Remove from this exterior product all possible pairs $\ell\ 
  \wedge r$ where $\ell$ is a leaf and $r$ is a root which are
  related by an edge.
\item Add to this exterior product pairs $ e ^+ \wedge e ^-$ for all
  edges $e$ between two inner vertices. Here $e ^+$ (resp. $e ^-$)
  stands for the upper (resp. lower) half-edge.
\end{enumerate}

The result is an exterior product on all half-edges of $F$ and an
auxiliary element $R_F$. One can assume that half-edges are gathered
by three according to their incident inner vertex. Replacing each such
triple $e_v^1\wedge e_v^2 \wedge e_v^3$ by the vertex $v$, one gets a
global orientation of $F$. One has to keep track of what has been
replaced. This is done by assigning the local orientation
$o_v=e_v^1\wedge e_v^2 \wedge e_v^3$ to the inner vertex $v$.

\smallskip

\begin{figure} \begin{center} \leavevmode
    \epsfig{file=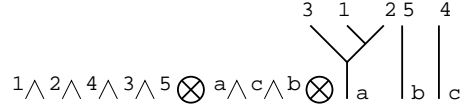,width=6cm} \caption{An outer-oriented
      forest on $\{1,2,3,4,5\}$.}\label{orfor}
  \end{center}\end{figure}

\begin{figure}  \begin{center}
    \leavevmode 
    \epsfig{file=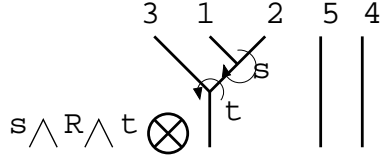,width=5cm} 
    \caption{An 
      inner-oriented forest on $\{1,2,3,4,5\}$.}\label{innexample}
  \end{center}\end{figure}

Here is an example of this equivalence of orientations. Consider the
outer-oriented forest shown in Fig. \ref{orfor}. One can compute the
corresponding inner orientation.
\begin{align*}
  & 1 \wedge 2 \wedge 4 \wedge 3 \wedge 5 \wedge R_F \wedge a \wedge c
  \wedge b 
\\ =\quad & 1 \wedge 2 \wedge 3 \wedge 5 \wedge R_F \wedge a \wedge b 
\\ =\quad & 1 \wedge 2  \wedge 3 \wedge R_F \wedge a
\\=\quad & 1 \wedge 2 \wedge 3
  \wedge R_F \wedge a \wedge e ^+ \wedge e ^- 
\\=\quad & (1 \wedge 2 \wedge e ^+)
  \wedge R_F \wedge (3 \wedge a \wedge e ^-),
\end{align*}
where $e ^+$ and $e ^-$ are the upper and lower half-edges of the
unique inner edge. Hence one can take the global orientation to be
$s\wedge R_F \wedge t$ (where $s$ is the upper vertex and $t$ the lower
one) and the local orientations to be $1 \wedge 2 \wedge e ^+$ at
vertex $s$ and $3 \wedge a \wedge e ^-$ at vertex $t$. The result is
shown in Fig. \ref{innexample}.

\smallskip

The grading is modified (but its parity is not changed) in order that
the forests with no inner vertex are in degree $0$, which will be
convenient in the next section. From now on, the degree of an
inner-oriented forest is the number of its inner vertices.

\subsection{Presentation of $\bess$}

From the known presentation of $B$, a presentation of $\bess$ by
generators and relations is given in this section.

Let $E_{i,j}$ be the inner-oriented forest with two trees on $\{i,j\}$
defined by the outer-oriented formula $E_{i,j}=(j\wedge i)\otimes
e_{i,j}$. It is symmetric of degree $0$. As an inner-oriented forest,
it is 
\begin{equation}
  R \otimes \stackrel{i}{|}\, \stackrel{j}{|}.
\end{equation}

\begin{figure}  \begin{center}
    \leavevmode \epsfig{file=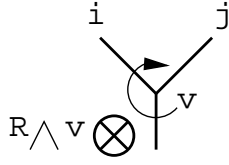,width=3cm}
    \caption{$\Omega_{i,j}$ as an inner-oriented tree.}\label{omega}
  \end{center}\end{figure}

Let $\Omega_{i,j}$ be the inner-oriented tree on $\{i,j\}$ defined by
the outer-oriented formula $\Omega_{i,j}=(i\wedge j)
\otimes\omega_{i,j}$. It is antisymmetric of degree $1$. As an
inner-oriented tree, it is given by Fig. \ref{omega}. 

\begin{prop}
  The operad $\bess$ is isomorphic to the quotient of the free operad
  on the generators $E_{i,j}$ symmetric of degree $0$ and
  $\Omega_{i,j}$ antisymmetric of degree $1$ by the relations
  \begin{align}
    \label{rel1_bess} E_{i,\star}\circ_\star E_{j, k}
    &=E_{k, \star }\circ_\star E_{i , j}, \\
    \label{rel2_bess}
    \Omega_{i,\star} \circ_\star E_{j,k} &= E_{j,\star} \circ_\star
    \Omega_{i,k} + E_{k,\star} \circ_\star \Omega_{i,j}.
  \end{align}
\end{prop}

\begin{proof}
  The tensor product by the operad $\Det$ acts essentially by changing
  all the signs. It is well known that the suspended operad has a
  presentation by similar generators and relations (up to sign) as it
  is simply given by a shift of grading at the level of algebras. Let
  us compute the new relations for our chosen generators. First,
  \begin{align*}
    E_{i,\star}\circ_\star E_{j, k }&=
    ((\star \wedge i)\otimes e_{i,\star})\circ_\star ((k\wedge j)\otimes
    e_{j,k})
    \\&= ((i \wedge \star)\circ_\star (k\wedge j))\otimes
    (e_{i,\star}\circ_\star e_{j,k})
    \\&= (i \wedge k \wedge j) \otimes (e_{i,\star}\circ_\star e_{j,k}).
  \end{align*}
  Therefore $E_{i,\star}\circ_\star E_{j, k}$ is invariant by cyclic
  permutation of $i,j,k$. One also has
  \begin{align*}
    \Omega_{i,\star} \circ_\star E_{j,k}&=((i\wedge \star)\otimes
    \omega_{i,\star})\circ_\star ((k\wedge j)\otimes
    e_{j,k})\\&=((i\wedge \star)\circ_\star(k\wedge
    j))\otimes(\omega_{i,\star}\circ_\star e_{j,k}) \\&=(i\wedge k
    \wedge j)\otimes (e_{j,\star}\circ_\star
    \omega_{i,k}-e_{k,\star}\circ_\star \omega_{i,j}).  \\&=(j\wedge i
    \wedge k)\otimes(e_{j,\star}\circ_\star \omega_{i,k}) +(k\wedge i
    \wedge j)\otimes(e_{k,\star}\circ_\star \omega_{i,j}))
    \\&=E_{j,\star}\circ_\star \Omega_{i,k}+E_{k,\star}\circ_\star
    \Omega_{i,j}.
  \end{align*}
\end{proof}

The action of $E$ is then described as follows. 
\begin{prop}
  Let $o_1 \otimes F_1$ and $o_2 \otimes F_2$ be two inner-oriented
  forests. Then
  \begin{equation}
    E_{\star,\#}\circ_\star (o_1 \otimes
    F_1) \circ_\# (o_2 \otimes F_2)= (o_1 \sqcup o_2) \otimes (F_1
    \sqcup F_2),
  \end{equation}
  where the global orientation $o_1 \sqcup o_2$ is obtained from
  $o_1\wedge r \wedge o_2$ by replacing $R_1\wedge r \wedge R_2$ by
  $R$. The local orientations are unchanged.
\end{prop}
\begin{proof}
  Let $o'_1\otimes o''_1$ and $o'_2\otimes o''_2$ be the corresponding
  outer orientations of $F_1$ and $F_2$. Using Prop. \ref{action_e},
  one has
  \begin{multline*}
    E_{\star,\#}\circ_\star (o_1 \otimes F_1) \circ_\# (o_2 \otimes
    F_2)=((\#\wedge \star) \otimes e_{\star,\#})\circ_\star
    (o'_1\otimes o''_1\otimes F_1) \circ_\# (o'_2\otimes o''_2\otimes
    F_2)\\=(-1)^{o'_1+o'_2+o'_2 o''_1}((\# \wedge \star )\circ_\star
   o'_1 \circ_\# o'_2)\otimes(e_{\star,\#}\circ_\star (
    o''_1\otimes F_1)\circ_\#(o''_2\otimes F_2))
    \\=(-1)^{(1+o''_1)(1+o'_2)}o'_1\wedge o'_2\otimes  o''_1\wedge o''_2
    \otimes (F_1 \sqcup F_2).
  \end{multline*}
  Hence the corresponding inner orientation is given by
  \begin{equation*}
    (-1)^{(1+o''_1)(1+o'_2)} o'_1\wedge o'_2 \wedge R \wedge o''_1
    \wedge o''_2.
  \end{equation*}
  
  On the other hand, let us compute the orientation corresponding to
  $o_1\sqcup o_2$.
  \begin{equation*}
    o'_1 \wedge R_1 \wedge o''_1 \wedge r \wedge o'_2 \wedge R_2
    \wedge o''_2 = (-1)^{(1+o''_1)(1+o'_2)} o'_1 \wedge o'_2 \wedge R
    \wedge o''_1 \wedge o''_2. 
  \end{equation*}
  Therefore the two orientations are the same.
\end{proof}

The action of $\Omega$ on trees has the following description.
\begin{prop}
  \label{omegatree}
  Let $o_1 \otimes T_1$ and $o_2 \otimes T_2$ be two inner-oriented
  trees. Then
  \begin{equation}
    \Omega_{\star,\#}\circ_\star (o_1 \otimes T_1)\circ_\#(o_2 \otimes
    T_2) = (o_1 \vee o_2) \otimes (T_1 \vee
    T_2),
  \end{equation}
  where the global orientation $o_1 \vee o_2$ is defined by
  $(-1)^{o_1} o_1 \wedge o_2$ modulo $R_1 \wedge R_2= R\wedge v$ where
  $v$ is the inner vertex of $\Omega$. The local orientations are
  unchanged.
\end{prop}
\begin{proof}
  Let $o'_1\otimes \rt_1$ and $o'_2\otimes \rt_2$ be the corresponding
  outer orientations of $T_1$ and $T_2$. Using Prop. \ref{action_o},
  one has
  \begin{multline*}
    \Omega_{\star,\#}\circ_\star (o'_1\otimes \rt_1 \otimes T_1)
    \circ_\# (o'_2\otimes \rt_2 \otimes T_2)\\=((\star \wedge
    \#)\circ_\star o'_1\circ_\#o'_2) \otimes
    (\omega_{\star,\#}\circ_\star \rt_1\otimes T_1 
    \circ_\# \rt_2\otimes T_2)\\=
    (-1)^{o'_1} (o'_1 \wedge o'_2) \otimes \rt \otimes (T_1 \vee T_2).
  \end{multline*}
  So the corresponding orientation is $(-1)^{o'_1} o'_1 \wedge o'_2
  \wedge R \wedge \rt$. Introducing pairs of half-edges gives
  \begin{equation*}
    (-1)^{o'_1} o'_1 \wedge o'_2 \wedge R \wedge \rt \wedge \rt_1
    \wedge e_1^- \wedge \rt_2 \wedge e_2^-,
  \end{equation*}
  where $e_1^-$ and $e_2^-$ are lower half-edges. This is equivalent
  with the local orientation $(\rt \wedge e_1^- \wedge e_2^-)$ at
  vertex $v$ (which is the local orientation of $\Omega$, see figure
  \ref{omega}) and orientation
  \begin{equation*}
    (-1)^{o'_1} o'_1 \wedge o'_2 \wedge R \wedge \rt_1 \wedge v \wedge \rt_2.
  \end{equation*}
  On the other hand, the proposed orientation is
  \begin{equation*}
    (-1)^{o_1} o'_1 \wedge R_1 \wedge \rt_1 \wedge o'_2 \wedge R_2
    \wedge \rt_2=(-1)^{o_1} o'_1 \wedge o'_2 \wedge R_1 \wedge \rt_1
    \wedge R_2
    \wedge \rt_2.
  \end{equation*}
  This matches the computed orientation, as $R_1 \wedge R_2 = R \wedge
  v$ and $(-1)^{o_1}=(-1)^{o'_1}$.
\end{proof}

\smallskip

Let us extend the definition of $\vee$ from trees to forests, as
follows. Let $F_1=T^1_1 \sqcup T^2_1 \sqcup \dots\sqcup T^m_1 $ and
$F_2=T^1_2 \sqcup T^2_2 \sqcup \dots \sqcup T^n_2$ be forests, where
the $T$ are trees. Define $F_1 \vee F_2$ to be the sum
\begin{equation*}
  \sum_{1 \leq a \leq m} \sum_{1 \leq b \leq n} (T^a_1 \vee T^b_2)
  \sqcup T_1^1 \sqcup \dots \sqcup\widehat{T}^a_1\sqcup \dots \sqcup
  T^2_2 \sqcup  \dots \sqcup\widehat{T}^b_2\sqcup \dots,
\end{equation*}
where $\widehat{T}$ means that this term is absent. In words, $F_1
\vee F_2$ is the sum over all possible pairings of a tree from $T_1$
and a tree from $T_2$, where these two trees are replaced in the
disjoint union $F_1 \sqcup F_2$ by their $\vee$ product.

\smallskip

Then Prop. \ref{omegatree} is still true for forests instead of just
trees, with the extended definition just given for $\vee$.

\begin{prop}
  Let $o_1 \otimes F_1$ and $o_2 \otimes F_2$ be two inner-oriented
  forests. Then
  \begin{equation}
    \Omega_{\star,\#}\circ_\star (o_1 \otimes F_1)\circ_\#(o_2 \otimes
    F_2) = (o_1 \vee o_2) \otimes (F_1 \vee
    F_2),
  \end{equation}
  where the global orientation $o_1 \vee o_2$ is defined by
  $(-1)^{o_1} o_1 \wedge o_2$ modulo $R_1 \wedge R_2= R\wedge v$ where
  $v$ is the inner vertex of $\Omega$. The local orientations are
  unchanged.
\end{prop}
\begin{proof}
  By recursion on the total number of trees in $F_1$ and $F_2$. The
  proposition is true if $F_1$ and $F_2$ are trees. Let us assume that
  $F_2$ has at least two trees. 

  One the one hand, 
  \begin{multline*}
    \Omega_{\star, \#}\circ_\star (o_1 \otimes F_1) \circ_\# ((o_2
    \sqcup o_3) \otimes (F_2 \sqcup F_3))
    \\=\Omega_{\star, \#}\circ_\star (o_1 \otimes F_1) \circ_\#
    (E_{\Delta,\infty} \circ_\Delta (o_2 \otimes F_2) \circ_{\infty}
    (o_3 \otimes F_3))
    \\= \Omega_{\star, \#} \circ_\# E_{\Delta,\infty} \circ_\star (o_1
    \otimes F_1) \circ_\Delta
    (o_2 \otimes F_2) \circ_{\infty} (o_3 \otimes F_3) 
    \\= (E_{\Delta,\#}\circ_\# \Omega_{\star,\infty} +
    E_{\infty,\#}\circ_\# \Omega_{\star, \Delta} ) \circ_\star (o_1
    \otimes F_1) \circ_\Delta (o_2
    \otimes F_2) \circ_{\infty} (o_3 \otimes F_3) 
    \\=(-1)^{o_2 o_3} E_{\Delta,\#}\circ_\# \Omega_{\star,\infty} \circ_\star
    (o_1 \otimes F_1)\circ_{\infty} (o_3 \otimes F_3) \circ_\Delta
    (o_2 \otimes F_2) \\+ E_{\infty,\#}\circ_\# \Omega_{\star, \Delta}
    \circ_\star (o_1 \otimes F_1) \circ_\Delta (o_2 \otimes F_2)
    \circ_{\infty} (o_3 \otimes F_3)
    \\= (-1)^{o_2 o_3} E_{\Delta,\#}\circ_\# ((o_1
    \vee o_3) \otimes (F_1 \vee F_3)) \circ_\Delta (o_2 \otimes F_2)
    \\+ E_{\infty,\#}\circ_\# ((o_1 \vee o_2) \otimes (F_1 \vee
    F_2)) \circ_{\infty} (o_3 \otimes F_3)
    \\=(-1)^{o_2 o_3} ((o_1 \vee o_3)\sqcup o_2) \otimes ((F_1 \vee
    F_3) \sqcup F_2) + ((o_1\vee o_2)\sqcup o_3) \otimes ((F_1 \vee
    F_2) \sqcup F_3).
  \end{multline*}

  On the other hand, the definition of $\vee$ implies that
  \begin{multline*}
    (o_1 \vee (o_2 \sqcup o_3))\otimes ( F_1 \vee (F_2 \sqcup F_3))
    \\= (o_1 \vee (o_2 \sqcup o_3))\otimes (( F_1 \vee F_3) \sqcup
    F_2)+(o_1 \vee (o_2 \sqcup o_3))\otimes ( (F_1 \vee F_2) \sqcup F_3)).
  \end{multline*}

  So it remains to compare the orientations. Using their defining
  properties, it is easy to see that 
  \begin{equation*}
    (-1)^{o_2 o_3} (o_1 \vee o_3)\sqcup o_2=o_1 \vee (o_2 \sqcup
    o_3)
    =(o_1\vee o_2)\sqcup o_3.
  \end{equation*}
  The proposition is proved.
\end{proof}

\section{A coproduct on $\bess$}
In this section, a map from $\bess$ to $\bess \otimes \bess$ is first
defined on generators, then shown to be given by an explicit formula.

\subsection{Definition on generators}
Let us define a coproduct $\Delta : \bess \to \bess\otimes \bess$ on
the generators $E_{i,j}$ and $\Omega_{i,j}$ of $\bess$ by
\begin{align}
  \label{coprod}
  \Delta(E_{i,j})&=E_{i,j} \otimes E_{i,j},\\
  \Delta(\Omega_{i,j})&=E_{i,j} \otimes \Omega_{i,j} +\Omega_{i,j}
  \otimes E_{i,j}.
\end{align}

\begin{prop}
  These formulas define a coassociative cocommutative morphism of
  operad from $\bess$ to $\bess \otimes \bess$, \textit{i.e.} the
  structure of a Hopf operad on $\bess$. In particular, each
  $\bess(I)$ inherits a structure of cocommutative coalgebra.
\end{prop}
\begin{proof}
  Coassociativity and cocommutativity are clear on generators. One has
  to check that the relations (\ref{rel1_bess}) and (\ref{rel2_bess})
  of $\bess$ are annihilated by $\Delta$. First,
  \begin{equation*}
    \Delta(E_{i,\star}\circ_\star E_{j,
      k})=(E_{i,\star}\otimes
    E_{i,\star}) \circ_\star ( E_{j, k} \otimes E_{j, k} )
    =(E_{i,\star}\circ_\star E_{j, k})\otimes (E_{i,\star}\circ_\star
    E_{j, k}), 
  \end{equation*}
  which inherits the invariance of $E_{i,\star}\circ_\star E_{j, k}$
  under cyclic permutations of $i,j,k$. Hence $\Delta$ vanishes on the
  relation (\ref{rel1_bess}).  For the other relation, on the one hand
  \begin{align*}
    \Delta&(\Omega_{i,\star} \circ_\star E_{j,k})= (E_{i,\star}
    \otimes \Omega_{i,\star} +\Omega_{i,\star} \otimes
    E_{i,\star})\circ_\star (E_{j,k} \otimes    E_{j,k})\\
    &=(E_{i,\star} \otimes \Omega_{i,\star}) \circ_\star (E_{j,k}
    \otimes E_{j,k})+(\Omega_{i,\star} \otimes
    E_{i,\star})\circ_\star (E_{j,k} \otimes    E_{j,k})\\
    &=(E_{i,\star}\circ_\star E_{j,k} )\otimes
    (\Omega_{i,\star}\circ_\star E_{j,k})+
    (\Omega_{i,\star}\circ_\star E_{j,k})\otimes
    (E_{i,\star}\circ_\star E_{j,k} )\\ &=(E_{i,\star}\circ_\star
    E_{j,k} ) \otimes (E_{j,\star} \circ_\star \Omega_{i,k})
    +(E_{i,\star}\circ_\star E_{j,k} ) \otimes (E_{k,\star}
    \circ_\star \Omega_{i,j})\\& +(E_{j,\star} \circ_\star
    \Omega_{i,k})\otimes (E_{i,\star}\circ_\star E_{j,k}
    )+(E_{k,\star} \circ_\star \Omega_{i,j}) \otimes
    (E_{i,\star}\circ_\star E_{j,k} ).
  \end{align*} 
  On the other hand,
  \begin{align*}
    \Delta&( E_{j,\star} \circ_\star \Omega_{i,k} ) =(
    E_{j,\star}\otimes E_{j,\star})\circ_\star
    (E_{i,k}\otimes\Omega_{i,k}+\Omega_{i,k} \otimes E_{i,k})\\
    &=(E_{j,\star} \circ_\star E_{i,k})
    \otimes(E_{j,\star}\circ_\star\Omega_{i,k} )
    +(E_{j,\star}\circ_\star\Omega_{i,k} ) \otimes (E_{j,\star}
    \circ_\star E_{i,k})\\
    &=(E_{i,\star} \circ_\star E_{j,k})
    \otimes(E_{j,\star}\circ_\star\Omega_{i,k} )
    +(E_{j,\star}\circ_\star\Omega_{i,k} ) \otimes (E_{i,\star}
    \circ_\star E_{j,k}),
  \end{align*}
  and a similar formula holds for $\Delta(E_{k,\star} \circ_\star
  \Omega_{i,j})$. From these formulas, it is clear that $\Delta$
  vanishes on relation (\ref{rel2_bess}).  This proves the
  proposition.
\end{proof}


\subsection{A poset on forests}

There is an explicit formula for the coproduct, which is a sum over
subsets of the set of inner vertices. A poset on forests involved in
this formula is described first.

A leaf is an \textit{ancestor} of a vertex if there is path from the
leaf to the root going through the vertex.

\smallskip

Let $F$ and $F'$ be two forests on the set $I$. Then $F' \leq F$ if
there is a topological map from $F'$ to $F$ with the following
properties :
\begin{enumerate}
\item It is increasing with respect to orientation towards the root.
\item It maps inner vertices to inner vertices injectively.
\item It restricts to the identity on leaves.
\end{enumerate}

In fact, such a topological map from $F'$ to $F$ is determined by the
image of inner vertices of $F'$. Indeed one can recover the map by
joining the image of an inner vertex with its ancestor leaves in $F'$.

\smallskip

\begin{figure}
  \begin{center}
    \leavevmode 
    \epsfig{file=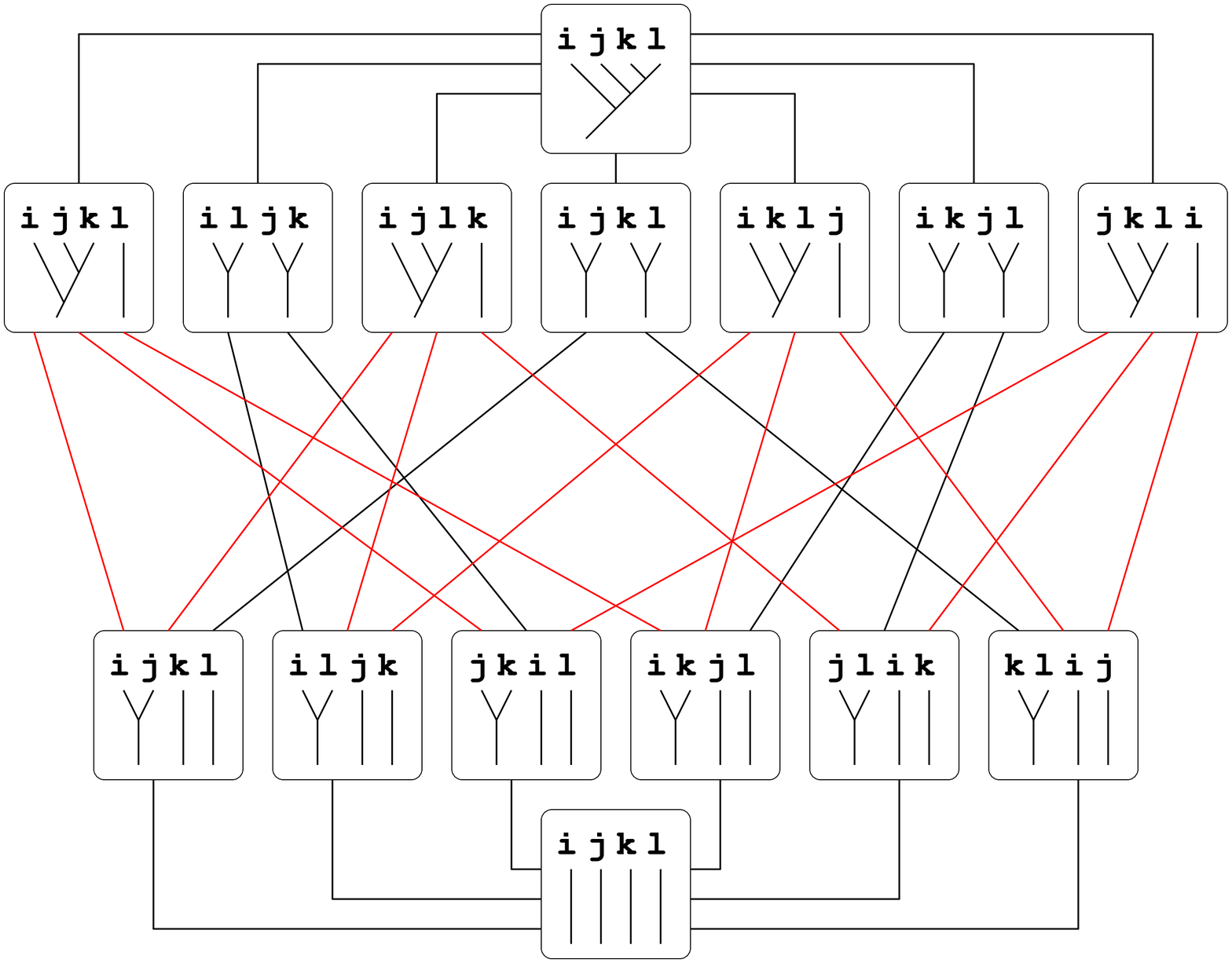,width=7cm} 
    \caption{An interval in the poset of forests on $\{i,j,k,\ell\}$.}
    \label{interval}
  \end{center}
\end{figure}

This relation defines a partial order on the set of forests on $I$.
The maximal elements of this poset are the trees. This poset is ranked
by the number of inner vertices. Fig. \ref{interval} displays an interval
in the poset of forests on the set $\{i,j,k,\ell\}$.

\smallskip

Remark : As can be seen on Fig. \ref{interval}, the interval in this
poset between the minimal element and one of the ``comb'' trees (which
have a leaf with all the inner vertices belonging to its path to the
root) can be identified to the partition lattice. The proof is by
identifying a forest with the partition of the set of leaves defined
by its trees. Details will be given elsewhere.

\smallskip

If $F$ is a forest on the set $I$ and $V$ is a subset of the set
$V(F)$ of inner vertices of $F$, let $\gamma(F,V)$ be the sum of
forests $F'$ such that $F' \leq F$ and the inner vertices of $F'$ are
identified with the elements of $V$. The sum $\gamma(V,F)$, which is
an element of the free $\ZZ$-module generated by the set of forests on
the finite set $I$, can also be considered as a set, as it has no
multiplicity. Indeed, there is at most one way to complete a injection
of inner vertices into a topological map from a given forest $F'$ to a
given forest $F$.

\begin{lemma}
  \label{gammadouble}
  Let $\sqcup$ and $\vee$ be the bilinear extensions of the operations
  $\sqcup$ and $\vee$ on forests.
  \begin{enumerate}
  \item Let $T=T_1 \vee T_2$ be a tree and $V'=\{v\} \sqcup V'_1
    \sqcup V'_2$ be a subset of $V(T)$ containing the bottom vertex $v$.
    Then $\gamma(T,V')=\gamma(T_1,V'_1)\vee \gamma(T_2,V'_2)$.
  \item Let $T=T_1 \vee T_2$ be a tree and $V'= V'_1 \sqcup V'_2$ be a
    subset of $V(T)$ not containing the bottom vertex $v$. Then
    $\gamma(T,V')=\gamma(T_1,V'_1)\sqcup \gamma(T_2,V'_2)$.
  \item Let $F=F_1 \sqcup F_2$ be a forest and $V'= V'_1 \sqcup V'_2$
    be a subset of $V(F)$. Then $\gamma(F,V')=\gamma(F_1,V'_1)\sqcup
    \gamma(F_2,V'_2)$.
  \end{enumerate}
\end{lemma}
\begin{proof}
  The second and third cases are essentially the same and easy
  consequences of the definition of the poset. If two sets $V_1$,
  $V_2$ of inner vertices of a forest $F$ have no ancestor leaf in
  common, then the set $\gamma(F,V_1 \sqcup V_2)$ is in bijection with
  the product $\gamma(F,V_1)\times \gamma(F,V_2)$. For $\gamma$ seen
  as a sum, this gives the expected result.
  
  The first case now. Any element of $\gamma(T,V')$ is a forest $F$
  with inner vertices $V'$. This forest can be restricted to $V'_1$
  and to $V'_2$ to give two forests $F_1$ and $F_2$. To be able to
  recover the forest $F$ from $F_1$ and $F_1$, it is necessary and
  sufficient to know to which tree of $F_1$ and to which tree of $F_2$
  the vertex $v$ was connected in $F$. Therefore the set
  $\gamma(T,V')$ is in bijection with the set of quadruples
  $(F_1,\alpha,F_2,\beta)$ where $F_1$ and $F_2$ are in
  $\gamma(T_1,V'_1)$ and $\gamma(T_2,V'_2)$, $\alpha$ is a tree of
  $F_1$ and $\beta$ is a tree of $F_2$.
  
  Therefore, seen as a sum, $\gamma(T,V')$ is exactly given by the
  bilinear extension of the operation $\vee$ on forests, which is a
  sum over the set of pairs of subtrees.
\end{proof}

\subsection{Explicit formula for the coproduct}

\begin{prop}
  \label{prop-coproduit}
  Let $o\otimes F$ be an inner-oriented forest. Then
  \begin{equation}
    \label{coproduit}
    \Delta(o \otimes F)=\sum_{V(F)=V' \sqcup V''}
    \big(o' \otimes \gamma(F,V' )\big)\otimes
    \big(o''\otimes \gamma(F,V'')\big).
  \end{equation}
  where the local orientations are unchanged and the global
  orientations satisfy $o' \wedge r \wedge o''=o$ modulo $R' \wedge r
  \wedge R''=R$.
\end{prop}
\begin{proof}
  The proof is a recursion on the number of inner vertices. The
  proposition is clear for trees with no inner vertex. The proof of the
  recursion step is done separately for trees and for forests with at
  least two trees.

  \paragraph{The case of trees}
  
  Let $o_1\otimes T_1$ and $o_2\otimes T_2$ be two inner-oriented
  trees and let $o\otimes T=(o_1 \vee o_2)\otimes (T_1 \vee T_2)$.
  Then
  \begin{multline}
    \label{deltatree}
    \Delta(o\otimes T)= \Delta\big(\Omega_{\star,\#} \circ_\star
    (o_1\otimes T_1) \circ_\# (o_2\otimes T_2)\big)
    \\ =\sum_{V(T_1)=V'_1
      \sqcup V''_1}\sum_{V(T_2)=V'_2 \sqcup V''_2} (\Omega_{\star,\#}\otimes
    E_{\star,\#}+E_{\star,\#}\otimes \Omega_{\star,\#}) \\ \circ_\star
    (o'_1 \otimes\gamma'_1 \otimes o''_1
    \otimes\gamma''_1)  \circ_\# (o'_2
    \otimes\gamma'_2 \otimes o''_2 \otimes\gamma''_2),
  \end{multline}
  where $\gamma ^*_i$ stands for $\gamma(T_i,V ^*_i )$.
  
  The first half of this formula corresponding to the
  expansion of the composition in $\Omega_{\star,\#}\otimes
  E_{\star,\#}$ is given by
  \begin{multline}
    \label{firsthalf1}
    \sum_{V(T_1)=V'_1 \sqcup V''_1}\sum_{V(T_2)=V'_2 \sqcup
      V''_2}(-1)^{o'_2 o''_1} \\(\Omega_{\star,\#}\circ_\star(o'_1
    \otimes\gamma'_1)\circ_\#(o'_2 \otimes\gamma'_2))\otimes
    (E_{\star,\#} \circ_\star(o''_1 \otimes\gamma''_1)\circ_\#(o''_2
    \otimes\gamma''_2))\\=\sum_{V(T_1)=V'_1 \sqcup
      V''_1}\sum_{V(T_2)=V'_2 \sqcup V''_2}(-1)^{o'_2 o''_1}
    \bar{o}' \otimes (\gamma'_1 \vee
    \gamma'_2)\otimes \bar{o}'' \otimes
    (\gamma''_1 \sqcup \gamma''_2),
  \end{multline} 
  the orientations satisfying
  \begin{alignat*}{2}
    o'_1 \wedge r_1 \wedge o''_1 &= o_1 \qquad &
    R'_1 \wedge r_1 \wedge R''_1 &= R_1 \\
    o'_2 \wedge r_2 \wedge o''_2 &= o_2 \qquad &
    R'_2 \wedge r_2 \wedge R''_2 &= R_2 \\
    (-1)^{o'_1} o'_1  \wedge o'_2 &= \bar{o}' \qquad &
    R'\wedge s & = R'_1 \wedge R'_2 \\
    o''_1 \wedge r''\wedge o''_2&=\bar{o}''\qquad &
    R''_1 \wedge r'' \wedge R''_2&=R''.
  \end{alignat*}
  
  On the other hand, one has to compute
  \begin{equation}
    \label{firsthalf2}
     \sum_{{V(T)=V' \sqcup V''} \atop {v \in V'}} o'\otimes
    \gamma(T,V') \otimes o''\otimes\gamma(T,V'').
  \end{equation}
  As $V(T)=\{v\} \sqcup V(T_1) \sqcup V(T_2)$, one can replace the sum
  by a double sum, using Lemma \ref{gammadouble} :
  \begin{equation*}
   \sum_{V(T_1)=V'_1
      \sqcup V''_1}\sum_{V(T_2)=V'_2 \sqcup V''_2} o'\otimes
    (\gamma'_1 \vee \gamma'_2 ) \otimes o'' \otimes (\gamma''_1 \sqcup
    \gamma''_2),
  \end{equation*}
  with the orientations determined by
  \begin{alignat*}{2}
    o'\wedge r \wedge o''&= o_1 \vee o_2 \qquad &
    R'\wedge r \wedge R''&=R \\
    (-1)^{o_1} o_1 \wedge o_2 &= o_1 \vee o_2\qquad &
    R_1 \wedge R_2 &= R \wedge s.
  \end{alignat*}
  
  All these conditions on orientations together imply that the
  orientations $o'\otimes o''$ and $(-1)^{o'_2 o''_1} \bar{o'} \otimes
  \bar{o}''$ are the same. Therefore (\ref{firsthalf1}) and
  (\ref{firsthalf2}) are equal.
  
  The other half of the sum (\ref{deltatree}), corresponding to the
  expansion of the composition in $E_{\star,\#}\otimes
  \Omega_{\star,\#}$, is shown in the same way to be equal to
  \begin{equation}
    \sum_{V(T)=V' \sqcup V''\atop {v \in V''}} o'\otimes \gamma(T,V') \otimes
    o''\otimes\gamma(T,V'').
  \end{equation}
  
  Therefore the full sum (\ref{deltatree}) is given by the expected
  formula (\ref{coproduit}) and the recursion step is done for trees.

  \paragraph{The case of true forests}
  
  Let $o_1\otimes F_1$ and $o_2\otimes F_2$ be two inner-oriented
  forests and let $o\otimes F=(o_1 \sqcup o_2)\otimes (F_1 \sqcup
  F_2)$. One has
  \begin{multline}
    \label{forestrecu1}
     \Delta(o \otimes F)=\Delta(E_{\star,\#} \circ_\star o_1\otimes F_1
    \circ_\# o_2\otimes F_2)\\ = \sum_{V_1=V'_1 \sqcup
      V''_1}\sum_{V_2=V'_2 \sqcup V''_2} (E_{\star,\#}\otimes
    E_{\star,\#}) \circ_\star (o'_1\otimes \gamma'_1 \otimes o''_1
    \otimes \gamma''_1) \circ_\# (o'_2\otimes \gamma'_2 \otimes o''_2
    \otimes \gamma''_2) \\ = \sum_{V_1=V'_1 \sqcup V''_1}\sum_{V_2=V'_2
      \sqcup V''_2} (-1)^{o'_2 o''_1} ( E_{\star,\#} \circ_\star
    o'_1\otimes \gamma'_1 \circ_\#o'_2\otimes \gamma'_2 )\otimes (
    E_{\star,\#} \circ_\star o''_1 \otimes \gamma''_1\circ_\# o''_2
    \otimes \gamma''_2) \\ = \sum_{V_1=V'_1 \sqcup V''_1}\sum_{V_2=V'_2
      \sqcup V''_2}(-1)^{o'_2 o''_1} \bar{o}'\otimes(\gamma'_1 \sqcup
    \gamma'_2) \otimes \bar{o}''\otimes(\gamma''_1 \sqcup \gamma''_2),
  \end{multline}
  where $\gamma ^*_i$ stands for $\gamma(F_i,V ^*_i )$ and the
  orientations satisfy
  \begin{alignat*}{2}
    o'_1 \wedge r_1 \wedge o''_1&=o_1 \qquad &
    R'_1 \wedge r_1 \wedge R''_1&=R_1\\
    o'_2 \wedge r_2 \wedge o''_2&=o_2\qquad &
    R'_2 \wedge r_2 \wedge R''_2&=R_2\\
    o'_1 \wedge r' \wedge o'_2&=\bar{o}'\qquad &
    R'_1 \wedge r' \wedge R'_2&=R'\\
    o''_1 \wedge r''\wedge o''_2&=\bar{o}''\qquad &
    R''_1 \wedge r'' \wedge R''_2&=R''.
  \end{alignat*}
  On the other hand, one has to compute
  \begin{equation}
    \label{forestrecu2}
    \sum_{V(F)=V'\sqcup V''} o'\otimes \gamma(F,V')
    \otimes o''\otimes \gamma(F,V'').
  \end{equation}
  As $V(F)=V(F_1) \sqcup V(F_2)$, one can replace the summation by two
  separate summations, using Lemma \ref{gammadouble} :
  \begin{equation}
    \sum_{V(F_1)=V'_1 \sqcup  V''_1}\sum_{V(F_2)=V'_2 \sqcup V''_2}
    o' \otimes (\gamma'_1 \sqcup \gamma'_2) \otimes o''\otimes
    (\gamma''_1 \sqcup  \gamma''_2),
  \end{equation}
  with the orientations satisfying
  \begin{alignat*}{2}
    o'\wedge r \wedge o''&=o_1 \sqcup o_2 & \qquad
    R'\wedge r \wedge R''&= R\\
    o_1 \wedge r_{12} \wedge o_2  &= o_1 \sqcup o_2 & \qquad
    R_1 \wedge r_{12} \wedge R_2 &= R.
  \end{alignat*}
  One can then show by using all the conditions above that the
  orientations $o'\otimes o''$ and $(-1)^{o'_2 o''_1} \bar{o}' \otimes
  \bar{o}''$ are the same, which implies that (\ref{forestrecu1}) and
  (\ref{forestrecu2}) are equal. The recursion step is done for
  forests.
  
  The proposition is proved.
\end{proof}

\begin{prop}
  The projection to the one-dimensional degree zero component is
  a counit. The inclusion of this degree zero component is
  an augmentation.
\end{prop}
\begin{proof}
  For a finite set $I$, there is just one forest of degree zero, which
  has no inner vertex. By inspection of the formula for the coproduct,
  this forest is grouplike. The second part of the proposition
  follows. This forest can only be obtained in the coproduct of $F$
  for the two summands given by $V(F) \sqcup \emptyset$ and $\emptyset
  \sqcup V(F)$, and the counit property is easily checked.
\end{proof}


\section{Algebras of labeled binary trees}

As it is sometimes more convenient to work with algebras rather than
coalgebras, we introduce here the algebra structure on the dual vector
space of $\bess(I)$.

\subsection{Description and properties}

Let us consider the dual basis, still indexed by inner-oriented
forests, of the dual vector space $\bess ^*(I)$, defined by the
following pairing from $\bess(I)\otimes\bess ^*(I)$ to $\QQ$.
\begin{equation*}
  \langle o\otimes F , o'\otimes F' \rangle =
  \begin{cases}
    0 \text{ if }F\not= F',\\
    1 \text{ if }F=F'\text{ and }o=o'.
  \end{cases}
\end{equation*}

The induced pairing from $\bess(I)\otimes\bess(I)\otimes\bess
^*(I)\otimes\bess ^*(I)$ to $\QQ$ is denoted again by $\langle \,
\rangle$.

\smallskip

It appears to be more convenient to use the opposite of the dual product.

\begin{prop}
  The opposite of the dual product is given by
  \begin{equation}
    (o_1\otimes F_1)
    \times (o_2\otimes F_2)= \sum_{(F,V_1 \sqcup V_2)} o\otimes F,
  \end{equation}
  where the orientations satisfy $o_1 \wedge r \wedge o_2 = o$ and
  $R_1 \wedge r \wedge R_2 = R$, the sum being over the set of pairs
  $(F,V_1\sqcup V_2)$ where $F$ is a forest and $V(F)=V_1 \sqcup V_2$
  a partition of the set of inner vertices of $F$ such that $F_1$
  appears in $\gamma(F,V_1)$ and $F_2$ appears in $\gamma(F,V_2)$.
\end{prop}
\begin{proof}
  The defining property of the dual product $\times^{op}$ is
  \begin{equation}
    \label{dualiser}
    \langle o_1 \otimes F_1 \otimes o_2 \otimes F_2 , \Delta( o
        \otimes F) \rangle
        =\langle ( o_1 \otimes F_1) \times^{op} (o_2 \otimes F_2) , o
        \otimes F \rangle.
  \end{equation}
  Let $ \Delta(o \otimes F)=\sum o'\otimes \gamma'\otimes o''\otimes
  \gamma''$, with the orientations given by $o' \wedge r \wedge o''=o$
  and $R'\wedge r \wedge R''=R$. The left hand-side of
  (\ref{dualiser}) can be computed as follows.
  \begin{multline*}
    \sum \langle o_1 \otimes F_1 \otimes o_2 \otimes F_2 , o'\otimes
    \gamma'\otimes o''\otimes \gamma'' \rangle\\ = \sum (-1)^{o' o_2}
    \langle o_1 \otimes F_1 ,o'\otimes \gamma' \rangle \langle o_2
    \otimes F_2,o''\otimes \gamma'' \rangle \\ =(-1)^{o_1
      o_2}\delta_{o',o_1} \delta_{ F_1 \in \gamma'} \delta_{o'',o_2}
    \delta_{ F_2 \in \gamma''},
  \end{multline*}
  where $\delta_{ F \in \gamma}$ is $1$ if $F$ belongs to the set/sum
  $\gamma$ and else $0$, and the orientations are identified in an
  obvious way. Here was used the fact that the sum $\gamma(F,V)$ is
  without multiplicity.
  
  Therefore, as taking the opposite product exactly removes the sign
  $(-1)^{o_1 o_2}$, one has
  \begin{equation*}
    \langle ( o_1 \otimes F_1) \times (o_2 \otimes F_2) , o
        \otimes F \rangle = \delta_{o',o_1} \delta_{ F_1 \in \gamma' }
        \delta_{o'',o_2} \delta_{ F_2 \in \gamma''}.
  \end{equation*}
  The proposition follows.
\end{proof}

Let the \textit{support} of a forest $F$, denoted by $\supp(F)$, be
the set of leaves which are not linked to the root by an edge,
\textit{i.e.} such that the path to the root contains at least one
inner vertex.

\begin{lemma}
  \label{support}
  Let $F$ be any forest appearing in the product of $F_1$ and
  $F_2$. Then $\supp(F)=\supp(F_1)\cup \supp(F_2)$.
\end{lemma}

\begin{prop}
  Let $F_1$, $F_2$ be two forests with disjoint supports. Then the only
  forest appearing in the product of $F_1$ and $F_2$ is the forest $F$
  with $\supp(F)=\supp(F_1)\sqcup \supp(F_2)$ which coincides with
  $F_1$ and $F_2$ on their respective support.
\end{prop}
\begin{proof}
  Any forest appearing in the product should have support the disjoint
  union of supports. The condition that $F_1 \leq F$ implies that the
  number of inner vertices of $F$ which are linked to the support of
  $F_1$ is greater or equal than the number of inner vertices of
  $F_1$. The same is true for $F_2$ and its support. But the number of
  inner vertices of $F$ is the sum of those of $F_1$ and $F_2$,
  therefore there is equality and the proposition follows.
\end{proof}

Let $Y_{i,j}$ be the element of $\bess^*(I)$ corresponding to the
forest with one inner vertex, with support $\{i,j\}$ and orientation
as in Fig. \ref{omega}.

\begin{lemma}
  \label{path}
  Let $F$ be a forest with $j\not\in \supp(F)$. Then the forests
  appearing in $F \times Y_{i,j}$ are exactly all forests
  obtained from $F$ by grafting a leaf $j$ to any edge in the path
  from $i$ to the root.
\end{lemma}
\begin{proof}
  It is clear that each such forest do appear in the product. We need
  only to show that there are no others. The forests which appear
  should have a vertex with leaves $i$ and $j$ as ancestors. As $j$ do
  not belong to the support of $F$, this vertex should be added to
  $F$. It can only be added on an edge of the path from $i$ to the
  root.
\end{proof}

\subsection{Some relations and open questions}

\begin{figure}
  \begin{center}
    \leavevmode 
    \epsfig{file=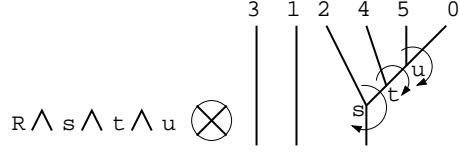,width=6cm} 
    \caption{$\LLL(2,4,5,0)$ on $\{0,1,2,3,4,5\}.$}
    \label{defll}
  \end{center}
\end{figure}

Let us introduce some notation. Let $\LLL(i,j,k,\ell)$ be the
inner-oriented forest on any set $I$ containing $\{i,j,k,\ell\}$,
which is defined on its support $\{i,j,k,\ell\}$ by the same
orientations and tree as Fig. \ref{defll}.

\begin{figure}
  \begin{center}
    \leavevmode 
    \epsfig{file=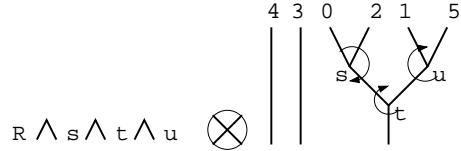,width=6cm} 
    \caption{$\YY(0,2,1,5)$ on $\{0,1,2,3,4,5\}.$}
    \label{defyy}
  \end{center}
\end{figure}

Let $\YY(i,j,k,\ell)$ be the inner-oriented forest on any set $I$
containing $\{i,j,k,\ell\}$, which is defined on its support
$\{i,j,k,\ell\}$ by the same orientations and tree as Fig.
\ref{defyy}.

\begin{lemma}
\label{symmetry}
  One has
  \begin{align}
    \LLL(i,j,k,\ell)&=-\LLL(i,j,\ell,k),\\
    \YY(i,j,k,\ell)&=-\YY(i,j,\ell,k),\\
    \YY(i,j,k,\ell)&=\YY(k,\ell,i,j).
  \end{align}
\end{lemma}

\smallskip

The following relations are satisfied in any algebra $\bess^*(I)$.

\begin{prop}
  Let $i,j,k$ be three distinct elements of $I$. Then
  \begin{equation}
    Y_{i,j}\times Y_{j,k}\times Y_{k,i}=0.
  \end{equation}
\end{prop}
\begin{proof}
  There can be no forest with $3$ inner vertices and support of
  cardinal $3$. The proposition therefore follows from lemma
  \ref{support}.
\end{proof}

\smallskip

Remark that $Y_{i,j}\times Y_{j,k}\times
Y_{k,\ell}=Y_{\ell,k}\times Y_{k,j}\times
Y_{j,i}$.

\begin{prop}
  Let $i,j,k,\ell$ be four distinct elements of $I$. Then
  \begin{equation}
    \sum Y_{i_1,i_2}\times Y_{i_2,i_3}\times Y_{i_3,i_4}=0,
  \end{equation}
  where the sum is over the set of total orders on $\{i,j,k,\ell\}$ up
  to reversal.
\end{prop}
\begin{proof}
  Using the product rule for the orientations and Lemma \ref{path}, one
  computes
  \begin{align*}
    Y_{i,j}\times Y_{j,k}\times Y_{k,\ell}&=
    \LLL(i,j,k,\ell)+\LLL(i,\ell,j,k)\\&+\LLL(\ell,i,k,j)
    +\LLL(\ell,k,j,i)+\YY(i,j,k,\ell).
  \end{align*}
  The sum of all $12$ similar terms obtained from this one by
  permutations of $\{i,j,k,l\}$ is then seen to vanish, using the
  antisymmetry and symmetry properties of $\LLL$ and $\YY$ stated
  in Lemma \ref{symmetry}.
\end{proof}

\smallskip

It is an interesting open problem to give a presentation by generators and
relations of the algebras $\bess ^*(I)$. 

\begin{question}
  Do the elements $Y_{i,j}$ generate $\bess ^*(I)$ ?
\end{question}

Assuming an affirmative answer, one can then ask

\begin{question}
  Do the relations above give a presentation of $\bess ^*(I)$ ?
\end{question}

\subsection{Differential forms and hyperplane arrangement}

Let $I$ be a finite set and $\CC^I$ be the vector space with
coordinates $(x_i)_{i \in I}$. Let $\textsf{H}_I$ be union of all
hyperplanes $x_i-x_j=0$ for $i\not= j$ in the subspace $\sum_{i \in I}
x_i=0$ of $\CC^I$.

\smallskip

It is well known from the work of Cohen (see \cite{cohen76,voronov})
that the Gerstenhaber operad is the homology of the little discs
operad, whose underlying spaces are homotopy equivalent to the
complement of the complex arrangements $\textsf{H}_I$. Therefore, by
the classical theorem of Arnold \cite{arnold} computing the cohomology
of this complement of arrangement, the coalgebra associated to a
finite set $I$ defined by the Hopf structure of the Gerstenhaber
operad has the following description : it is isomorphic to the dual of
the subalgebra generated by all forms $d(x_i-x_j)/(x_i-x_j)$ for
$i\not= j$ in $I$ inside the algebra of differential forms on the
complement of $\textsf{H}_I$,

\smallskip

The differential forms $Y_{i,j}=d(x_i-x_j)/(x_i-x_j)^2$ for
$i\not= j$ in $I$ are defined on the complement of $\textsf{H}_I$.
Obviously, they satisfy $Y_{i,j}=-Y_{j,i}$.

\smallskip

Let $i,j,k$ be three distinct elements of $I$. Then one has clearly
\begin{equation}
    Y_{i,j}\wedge Y_{j,k}\wedge Y_{k,i}=0.
\end{equation}

Further experimental evidence has been obtained showing that the
algebra on forests of binary trees considered in this article should
be isomorphic to a quotient of the subalgebra generated by the
$Y_{i,j}$ inside the algebra of differential forms on the complement
of $\textsf{H}_I$.

\bibliographystyle{plain}
\bibliography{bessel}

\begin{thebibliography}{1}

\bibitem{arnold}
V.~I. Arnold.
\newblock The cohomology ring of the group of dyed braids.
\newblock {\em Mat. Zametki}, 5:227--231, 1969.

\bibitem{cohen76}
F.~R. Cohen.
\newblock The homology of {$\mathcal{C}_{n+1}$}-spaces, {$n\ge0$}.
\newblock In {\em The homology of iterated loop spaces}, volume 533 of {\em
  Lecture Notes}, pages 207--351. Springer-Verlag, 1976.

\bibitem{ginzkapr}
V.~Ginzburg and M.~Kapranov.
\newblock Koszul duality for operads.
\newblock {\em Duke Math. J.}, 76(1):203--272, 1994.

\bibitem{grosswald}
E.~Grosswald.
\newblock {\em Bessel polynomials}.
\newblock Springer, Berlin, 1978.

\bibitem{krallfrink}
H.~L. Krall and O.~Frink.
\newblock A new class of orthogonal polynomials: {T}he {B}essel polynomials.
\newblock {\em Trans. Amer. Math. Soc.}, 65:100--115, 1949.

\bibitem{markl}
M.~Markl.
\newblock Distributive laws and {K}oszulness.
\newblock {\em Ann. Inst. Fourier (Grenoble)}, 46(2):307--323, 1996.

\bibitem{livreMSS}
M.~Markl, S.~Shnider, and J.~Stasheff.
\newblock {\em Operads in algebra, topology and physics}.
\newblock American Mathematical Society, Providence, RI, 2002.

\bibitem{may}
J.~P. May.
\newblock Definitions: operads, algebras and modules.
\newblock In {\em Operads: Proceedings of Renaissance Conferences (Hartford,
  CT/Luminy, 1995)}, pages 1--7, Providence, RI, 1997. Amer. Math. Soc.

\bibitem{voronov}
A.~A. Voronov.
\newblock Homotopy {G}erstenhaber algebras.
\newblock In {\em Conf\'erence Mosh\'e Flato 1999, Vol. II (Dijon)}, pages
  307--331. Kluwer Acad. Publ., Dordrecht, 2000.

\end{thebibliography}
\end{document}